%% file: HADAHA-Arnold.tex
\documentclass[12pt]{amsart}
\usepackage{amsmath,amsfonts,amssymb,amscd,color,epsfig}
\usepackage{latexsym}

\newcommand{\comment}[1]{}

\title [Affine Hecke Algebras via DAHA] 
{Affine Hecke Algebras via DAHA}
\author[Ivan Cherednik]{Ivan Cherednik $^\dag$} 

\address[I. Cherednik]{Department of Mathematics, UNC 
Chapel Hill, North Carolina 27599, USA\\
chered@math.unc.edu}

\thanks{\dag
\ Partially supported by NSF grant
DMS--1363138}

\input{macro}

\font\eightrm=cmr8

\begin{document}
\par
{\centering 
Dedicated to Masaki Kashiwara, a great master of harmonic
analysis,  \\
on the occasion of his 70th birthday
\medskip
\par} 
\maketitle




\comment{
A method is suggested for obtaining the Plancherel measure 
for Affine Hecke Algebras as a limit of 
integral-type formulas for inner products in the polynomial 
and related modules of Double Affine Hecke Algebras.
The analytic continuation necessary here is a 
generalization of ``picking up 
residues" due to Arthur, Heckman, Opdam and others, which
can be traced back to Hermann Weyl. Generally, it is a
finite sum of integrals over double affine 
residual subtori; a complete formula is presented for $A_1$
in the spherical case. 
}

{
This paper is based on the lecture delivered at the conference
``Algebraic Analysis and Representation Theory" in
honor of Masaki Kashiwara's 70th birthday, an author's
prior talk at MIT and his course at UNC Chapel Hill. 
It is aimed at obtaining
the Plancherel measure for the regular 
representation of 
Affine Hecke Algebras (AHA) as the limit $q\to 0$ of
the integral-type formulas for the DAHA inner products in
the polynomial and related modules. The usual integral
formulas generally serve only $\Re k>0$ (in the DAHA 
parameters $t=q^k$) and
must be analytically continued to negative $\Re k$, which 
is a $q$-generalization of ``picking up  
residues" due to Arthur, Heckman, Opdam and others
(can be traced back to Hermann Weyl). When this is done, 
we arrive at finite sums of integrals over {\em double affine 
residual subtori\,}, though the full procedure is known only 
in type $A$ by now. This is not related to the DAHA 
reducibility of the polynomial and similar DAHA modules. The 
formulas are nontrivial for any
$\Re k<0$, not only for {\em singular\,} $k<0$ resulting
in the DAHA reducibility. For singular $k$ in type $A$, 
they provide the decomposition of the {\em polynomial 
representation\,} in terms of the irreducible modules.
As we demonstrate, this is quite interesting even for $A_1$. 
\vfil

The decomposition of the regular AHA representation
in terms of (unitary) irreducible modules is
an important part of algebraic harmonic analysis,
involving deep geometric methods 
(Kazhdan- Lusztig and others). As an expected 
application, our approach would allow
to interpret formal degrees of AHA discrete series
via DAHA, without any geometry. Paper 
\textsuperscript{\cite{O11}}
do this within the AHA theory, but the DAHA level is
expected to be quite clarifying and more powerful (with
an additional parameter $q$).  

We mainly discuss the spherical case and
provide explicit analytic continuations only for $A_1$.
The key is that the mere uniqueness of the DAHA inner product 
fixes uniquely the $q$-generalization of the
corresponding Arthur-Heckman-Opdam formula,
including very interesting  $q$-counterparts of formal degrees.

Even in the spherical case, 
the procedure of analytic continuation to $\Re(k)<0$
is technically involved. There are no significant theoretical 
challenges here, but practical finding double affine residual 
subtori and their contributions to the DAHA inner products
is performed (partially)
only for $A_n$ at the moment. The passage to the
whole regular representation will presumably require
the technique of {\em hyperspinors\,}, which we outline a 
bit at the end of this paper. 

Importantly, there is no 
canonical AHA-type trace in the DAHA theory; instead, 
we have the theory of {\em DAHA coinvariants\,}
serving DAHA anti-involutions. There are of course
other aspects of DAHA harmonic analysis: the unitary dual,
calculating Fourier transforms of DAHA modules and so on,
where the regular representation of DAHA must be studied
(by analogy with the AHA theory). However, we focus on 
the spherical part of the regular AHA representation, which
becomes an {\em irreducible\,} module in the DAHA theory. 

Only basic references are provided in the paper; see there 
for further information. Also, the general AHA and DAHA theory
is quite compressed. Full details are provided for 
$A_1$; the generalization to $A_n$ follows the same lines.


\let\thefootnote\relax\footnote{\dag
\ Partially supported by NSF grant
DMS--1363138}
}

\mathversion{normal}



We  do not give general definitions (for arbitrary root systems) 
in this note. These definitions, including the basic
features of DAHA inner products, are (published and) sufficiently 
well known for $\Re(k)>0$; the main references are
\textsuperscript{\cite{C101,C5,CMa}}. The extension of the
corresponding {\em integral\,} formulas to $\Re(k)<0$ is the
aim of this work; we think that the case of
$A_1$ gives a clear direction. Importantly,
only relatively elementary  tools 
from ``$q$-calculus" are needed here. DAHA theory guarantees that
the inner products under consideration do have analytic and
meromorphic continuations, but does not provide explicit 
formulas (summations of integrals over generalized residual
subtori). This is actually similar to 
\textsuperscript{\cite{HO,O11}}.
The sections of the paper are:
\vskip 0.1cm  

{
1. On Fourier Analysis \ \  \ \ 
2. AHA--decomposition 

3. Shapovalov pairs \ \ \ \ \ \ \ \ 
4. Rational DAHA ($A_1$)

5. General DAHA ($A_1$) \ \ \ 
6. Analytic continuation 

7. P--adic limit\ \ \ 
\ \ \ \ \ \ \ \ \ \ \ \  8. Conclusion.

\vskip 0.2cm
}

\comment{
\noindent
\textcolor{red}{\sf  \ Warmest congratulations to
Kashiwara sensei (70!)}

\textcolor{red}
{{\sc  a great master of harmonic analysis} !!}
\vskip 0.3cm
\vfill
}

The following table sketches the basic ``levels" in
harmonic analysis on AHA vs. those in the corresponding 
DAHA theory. We think it explains our take on the AHA 
theory sufficiently well. The second column here is
technically adding an extra parameter $q$ to the theory.
Conceptually, the passage to DAHA provides important 
rigidity, which is of clear importance for 
the AHA Plancherel formula and related problems. The key
is that Fourier 
transform is essentially an involution in the DAHA theory, 
which is so different from AHA theory and 
Harish-Chandra theory.
\vskip 0.2cm

\begin{center}
\begin{tabular}{|c|c|}
\hline
\textcolor{blue}{\sf HA on AHA} & 
\textcolor{blue}{\sf HA on DAHA}\\ 
\hline
Unitary(spherical) dual   & {
Polynomial/induced modules}\\   
\hline
AHA Fourier transform &  {
\HH-- automorphism $Y\!\!\rightarrow 
\!\!X^{\!-\!1}$} \\
\hline
{
 Trace formulas, $L^2(\h)$} & 
{
Inner products as integrals}\\
\hline 
\end{tabular} 
\end{center}
\vskip 0.3cm

The three ``stages" in the first column are
common in harmonic analysis on symmetric
spaces and related/similar theories. One can ask
the same questions for DAHA. The first two 
stages are meaningful; DAHA provides an important
source of new infinite-dimensional 
\textcolor{red}{\it unitary\,}
theories, which are of great demand in analysis
and physics. However, our understanding is that
there is no {\it canonical\,} 
DAHA trace; accordingly, it is not clear
what the theory of $L^2(\HH)$ can be. Instead, we
have the analytic theory of DAHA involutions and 
coinvariants. The key for us is the interpretation
of the spherical part of the regular AHA representation
as DAHA polynomial representation, with very rich 
structures.  
  
{
One of possible applications of this program can be a
new approach to  
{\it formal degrees of AHA discrete series\,}
via DAHA.} Let us mention (at least) 
{
Kazhdan, Lusztig, Reeder, Shoji, Opdam, Ciubotaru, 
S.Kato} in this regard; see some references 
below.

\section{\textcolor{blue}{\sc On Fourier Analysis}}

\noindent
This section is mainly needed to put this work 
into perspective. 
The classical Fourier Transform, 
{\sf FT}$=\int e^{2\la x}\{\cdot\}\,dx,$ 
is naturally associated with the automorphism 
$x\rightarrow y=d/dx\rightarrow -x$ of the $d=1$ {\em 
Heisenberg algebra}.
This can be readily extended to any dimensions $d$ and to any 
root systems. 
Its spherical generalization is a famous Harish-Chandra 
transform, but this is beyond the Heisenberg algebras. 
Classically, the Fourier transform can be related to
{\tiny 
$\left(\begin{array}{cc}0& 1 \\
                       -1 & 0 \\
\end{array} \right)$} $\in SL_2$, though this kind of
interpretation seems a special feature of $d=1$ ($A_1$); 
see below.
Similarly, the {\em Weyl algebra\,} 
at $q\!=\!e^{\frac{2\pi\imath}{N}}$
can be used to study
$F_N\!=\!\sum_{j=0}^{N-1} q^{\la j}\{\cdot\}$, which can be
readily generalized to any $d$. Finding counterparts of 
Heisenberg/Weyl algebras
{\em directly\,} serving the Harish-Chandra transform and
its variants and generalizations is a natural question here.
DAHA essentially manages this. Using Lie groups here
is generally insufficient even for the classical one-dimensional
hypergeometric function.

\vskip 0.2cm

\textcolor{blue}{\sf Some famous challenges.}
\vskip 0.1cm

Problem 1. Extending Lie theory  
from spherical functions to  
\textcolor{red}{hypergeometric
functions} (in any ranks), the Gelfand Program. 
Here Kac-Moody algebras (conformal blocks, to be more
exact) and Lie super-groups can be used, 
but the problem still appeared beyond Lie theory.
\vskip 0.2cm

Problem 2. Can \textcolor{red} {\it Fourier transform\,} 
be interpreted as a reflection in the Weyl group (in any rank)? 
Unlikely so. Say, there are $3$ candidates (reflections)
for {\sf FT} in $SL_3$, but it can be expected {\it unique\,}
due to the key property of {\sf FT} in
any theories: they send   
polynomials to $\delta$-functions.

\vskip 0.2cm
Problem 3. A counterpart of 
{\sf FT}$(e^{-x^2})\!=\!\sqrt{\pi}e^{+\la^2}$ 
at roots of unity is the formula
$F_N(q^{j^2})=$
\textcolor{red}{\it \mathversion{normal} $\zeta$}
$\sqrt{N}q^{-\la^2}$ 
for
$\zeta\in\{0,1,\imath,1+\imath\}$.
The Weyl algebra gives $\sqrt{N}$ but
does not catch \textcolor{red}{\it\mathversion{normal} $\zeta$},
i.e. it provides only the absolute values
of the \textcolor{red} 
{\it Gauss sums\,}.
Can this be improved?

\vskip 0.2cm
\textcolor{blue} {\sf DAHA approach.}
Concerning \textcolor{red}{\it Problems 1-2\,}, 
the reproducing kernel of the DAHA-Fourier transform 
(its square is essentially one) is the generalized 
``global" difference hypergeometric function; any root 
systems were managed. \textcolor{red}{\it Problem 3\,} can 
be settled too (within the theory of DAHA-Gauss-Selberg sums). 
The question we (partially) address in this paper is, 
\textcolor{red}{\it what do these developments give for AHA?} 
\vskip 0.1cm

Before coming to this, let us quickly discuss 
{\it global hypergeometric functions}, which seem the
main application of DAHA with known and expected
applications well beyond harmonic analysis. In the $p$-adic
limit, which is $q\to 0$, they reduce to the polynomials;
classical hypergeometric functions are for
$q\to 1$. 

\vspace*{0.1cm}

\textcolor{blue}{\sf Global functions\  $\Phi_{q,t}(X,\La)
,\  q<1$.}\textsuperscript{\cite{C5}}
These functions are defined as reproducing kernels of 
the DAHA Fourier Transform; no difference equations are
used in this approach for their definition. 
To calculate them we use that
this transform sends Laurent polynomials in terms of $X=q^x$ 
times the Gaussian $q^{-x^2/2}$
to such polynomials times  $q^{+x^2/2}$, which
is similar to theory of Hankel
transform. This gives an explicit formula for the series
$\tilde{\Phi}_{q,t}(X,\La)\!\equal\!
\theta_R(X)\theta_R(\La)\Phi_{q,t}(X,\La)/\theta_R(t^{\rho})$
in terms of $P_\mu(X)P_\mu(\La)$ for Macdonald
polynomials $P_\mu$.

\vskip 0.1cm
This is general theory, for any 
(reduced, irreducible)
root systems $R$. The Laurent polynomials are in terms of
$X_\la=q^{(x,\la)}$ for $\la\in P$ (the weight lattice
for $R$), $(\cdot,\cdot)$ is
the standard $W$-invariant inner product,
$x^2=(x,x)$, $\theta_R$ is the usual 
theta-series associated with $R$, 
$\mu\in P$. This is from \textsuperscript{\cite{C5}}; see
also \textsuperscript{\cite{C101,ChW}}.
The series $\tilde{\Phi}_{q,t}(X,\La)$
is  absolutely convergent as $|q|<1$, 
$W$-invariant with respect to both, $X$ and $\La$, and, 
importantly, $X\!\leftrightarrow\! \La$-symmetric (as for the
Bessel functions). Here one must avoid the poles 
of the coefficients of Macdonald polynomials and 
zeros of $\theta_R(X)$; otherwise the convergence is 
really ``global". 
Such global functions are missing in
the (differential) Harish-Chandra theory.

\vskip 0.1cm  
Furthermore, let $X=q^x,\La=q^\la$. Assume that
$\la=w(\la_+)$ for dominant $\la_+$ such that 
$\Re(\la_+,\al_i)>0$, i.e. that $\la$ is generic. 
Then $\Phi_{q,t}(X,\La)$ under some explicit normalization
becomes an asymptotic series
$\Phi^{a\!s}_{q,t}(X,\La)= q^{-(x,\la_+)}
t^{(x+\la_+,\rho)}(1+\ldots)$ 
as $\Re(x,\al_i)\to +\infty$. 
\vskip 0.1cm

\textcolor{blue}
{\sf Harish-Chandra decomposition.} 
See I.Ch\textsuperscript{\cite{ChW,ChO1}}, 
J.Stokman\textsuperscript{\cite{Sto2}}. It is:
$$\Phi_{q,t}(X,\La)= \sum_{w\in W}
\si_{q,t}(w(\La))\,
\Phi_{q,t}^{a\!s}(X,w(\La))\ 
$$  
for the $q,t$-extension $\si_{q,t}(\La)$
of the Harish-Chandra $c$-function.
Note that (naturally)
 $\l_p\Phi_{q,t}(X,\La)=p(\La)\Phi_{q,t}(X,\La)$
for $p\in \C[X]^W$ and Macdonald-Ruijsenaars operators
$\l_p$ in type $A$
(and due to DAHA for any root systems). This 
relation is not used in the definition
of $\Phi_{q,t}$. The {\em recovery
formula\,} is more important here: $P_\la$ is proportional to
$\Phi_{q,t}(X,\La)$ for $\La\!=\!t^{\rho} q^\la$
(with an explict coefficient of proportionality).  

\vskip 0.2cm
Let us give the exact formulas for $A_1$.
For Rogers-Macdonald polynomials
$P_n(X)$, $\mu$\, provided below, and 
$\theta(X)\equal\sum_{m=-\infty}^\infty q^{mx+m^2/4}$, 
\begin{align*}
&\frac{\theta(X)\theta(\La)}{\theta(t^{1/2})}\Phi=
\tilde{\Phi}_{q,t}(X,\La)\equal\,
\sum_{n=0}^\infty\,q^{\frac{n^2}{4}}\,
t^{\frac{n}{2}}
\,\frac{P_n(X)P_n(\La)\, (\mu)_{\hbox{\sc\tiny ct}}}
{(P_n P_n\,\mu)_{\hbox{\sc\tiny ct}}},\ |q|<1.
\end{align*}

For $|X|\!<\!|t|^{\frac{1}{2}}|q|^{-\frac{1}{2}}$,
the Harish-Chandra formula reads:\  
$\tilde{\Phi}_{q,t}(X,\La)\!=$
\vskip -0.35cm
\begin{align*}
=&\,(\mu)_{\hbox{\sc\tiny ct}}\,\si(\La)\,\,\theta(X\La t^{-1/2})
\,\sum_{j=0}^\infty\,
(\frac{q}{t})^{\!{}^{j}}\,X^{2j}\prod_{s=1}^j\frac{(1-tq^{s-1})
(1-q^{s-1}t\La^{-2})}{(1-q^s)(1-q^s\La^{-2})}
\notag\\
+&\,(\mu)_{\hbox{\sc\tiny ct}}\,
\si(\La^{-1})\theta(X\La^{-1}t^{-1/2})\,
\sum_{j=0}^\infty
(\frac{q}{t})^{\!{}^{j}}X^{2j}\prod_{s=1}^j\frac{(1\!-\!tq^{s-1})
(1\!-\!q^{s\!-\!1}t\La^{2})}{(1\!-\!q^s)(1\!-\!q^s\La^{2})},
\notag
\end{align*}
where $\si(\La)=\prod_{j=0}^\infty
\frac{1-tq^j \La^2}{1-q^j \La^2}$
is the $q,t$-generalization of the Harish-Chandra $c$-function;
$(\cdot)_{\hbox{\sc\tiny ct}}$ is the constant term. See
\textsuperscript{\cite{ChO1}}. 
The sums here are nothing but
(special) Heine's {\it basic} hypergeometric functions. 
Letting here $t\to 0$ in type $A_n$,
the asymptotic expansions of the resulting 
{\em global $q$-Whittaker function\,} are essentially the 
Givental-Lee functions. This is an important connection
between the physics $B$-model (the usage of the {\em global\,}
function) and the $A$-model (the usage of its asymptotic 
expansions). We note that $\Phi$ is actually an entirely algebraic
object, uniquely determined by
its asymptotic behavior, including the 
\textcolor{blue} {\it walls} (resonances),
when $\Re(\al,\la)=0$ for
some roots $\al$ (the theory of resonances is still incomplete). 

\section{\textcolor{blue}{\sc AHA--decomposition}}\label{SEC:AHA}

Let $R\subset \R^n$  be a root system,\,
$Q\subset P$(the weight lattice),
$W=\lan\!s_\al\!\!\ran$ for $\al\in R$,
$\tilde{W}=W\lsmash Q\subset \hat{W}=
W\lsmash P=\tilde{W}\lsmash \Pi$, where $\Pi\!=\!P/Q$.
 
Then $\h\equal \lan\Pi,\, T_i(0\le i\le n)\ran /\{$homogeneous
Coxeter relations for $T_i$,\and 
$(T_i-t^{\frac{1}{2}})(T_i+t^{-\frac{1}{2}})=0\for  
1\le i\le n$ $\}$, where $\R$ will be the ring of coefficients,
including $q,t^{\pm 1/2}$. This is convenient
to avoid the complex conjugation in the scalar products 
(and for positivity).  

\vskip 0.3cm

We set $T_{\hw}\!=\!\pi T_{i_l}\!\cdots\! T_{i_1}$ for 
reduced decompositions
$\hw\!=\!\pi s_{i_l}\!\cdots\! s_{i_1}\!\in\! \hat{W},$
where $l\!=\!l(\hw)$ is the length of $\hw$. 
The canonical anti-involution, trace and scalar product are: 
\vskip 0.2cm

 {\mathversion{normal}
$T_{\hw}^\star\equal T_{\hw^{-1}},\, 
\lan T_{\hw}\ran =\de_{id,\hw},$}\ 
$\lan f,g\ran \equal \lan f^\star g\ran=
\sum_{\hw\in \hat{W}}\,c_{\hw}d_{\hw},\, $ 

\vskip 0.1cm
\noindent
where
$f\!=\!\sum c_{\hw} T_{\hw},\ g\!=\!
\sum d_{\hw} T_{\hw}\, \in\, L^2(\h)=\{f,\, c_{\hw}\in \R,\ 
\sum c_{\hw}^2<\infty \}$.

\vskip 0.2cm
According to Dixmier, \textcolor{blue}{\it 
 $\lan f,g\ran=
\int_{\pi\in\h^\vee} \hbox{Tr}(\pi(f^\star g))d\nu(\pi)$.}  
We omit here some analytic details concerning
the classes of functions. In the spherical
case (referred to as ``sph" later on), one takes
$f,g\in P_+\h P_+$, where 
$P_+\equal\sum_{w\in W} t^{\frac{l(w)}{2}}T_{w}.$ 
The measure reduces correspondingly. 

\vskip 0.2cm
Macdonald found an integral formula for 
 $\nu_{sph}(\pi),$ as $t>1$. Its extension to $0<t<1$ 
(due to ... Arthur, Heckman-Opdam, ...) by the analyticity
is sometimes called ``picking up residues"
\,\textsuperscript{\cite{CKK,HO,O11,OS}}. The final
formula (for any $t$) generally reads:
 
\noindent
\textcolor{blue}\it
{$$\int \{\cdot\}\, d\nu_{sph}^{an}(\pi)= 
\sum C_{s,S}\cdot\int_{s+iS}\{\cdot\} \,d\nu_{s,S},$$}

\noindent 
summed over (affine)
\textcolor{blue}{\sf residual subtori} $s+S.$
Residual \textcolor{blue}{\it points} (very interesting and
the most difficult to reach) correspond to
square integrable irreducible modules
(as their characters $\chi_\pi$ extend to $L^2(\h)$).

This formula involves deep algebraic geometry, the  
Kazhdan-Lusztig theory
\textsuperscript{\cite{KL1,Lu}}. 
In our approach via DAHA, this very formula expected
to be a reduction of the analytic continuation of the
DAHA inner product in the integral form, which requires
only $q$-calculus. The main claim is as follows.
It is in the spherical case and is a theorem for any
(reduced, irreducible) root systems, with an important
reservation that the explicit formula is known
by now only in type $A$ (unpublished).   

\vskip 0.3cm


\noindent
\textcolor{red}{\it The $q,t$-generalization of the 
picking up residues is the presentation of  the 
inner product in the DAHA polynomial representation as
sum of integrals over DAHA residual subtori. 
Only the whole sum satisfies the DAHA invariance, and
the corresponding $C$-coefficients are uniquely determined 
by this property. Upon the limit $q\to 0$, this approach
potentially provides explicit formulas
for the $C_{s,C}$-coefficients above, including formal degrees
(for the residual points).}


\section{
\textcolor{blue}{\sc Shapovalov pairs}
\textsuperscript{\cite{CMa}}}
 
We will now switch to the DAHA harmonic analysis.
In contrast to the Harish-Chandra theory, where we
mainly have two theories based on the imaginary 
and real integration, the so-called {\em compact\,} and 
{\em non-compact\,}
cases, here we have more options. Let us try to outline them,
disregarding
various (many) specializations and 
the open project aimed at the passage from the $q$-Gamma
function in DAHA theory to the $p$-adic Gamma (this is
doable, but there are no works on this so
far). 

We think that there are
essentially $6$ major theories by now, corresponding
to different choices of ``integrations"; some connections
are shown by arrows. We stick to the imaginary integration
in this paper. 

\vskip 0.2cm
\centerline{
DAHA INTEGRATIONS:}

\begin{center}
\noindent
\begin{tabular}{|c||c|}
\hline
 imaginary ($|q|\neq 1)$            & real ($|q|\neq 1$)\\ 
 $\Downarrow$                 & $\Downarrow$  \\
 constant term ($\forall q$)  & Jackson sums \\ 
$\Uparrow$                    & $\Downarrow$  \\
 the case $|q|=1$ \hfill  $\Rightarrow$ 
                              & $\Rightarrow$ \hfill
                                 roots of unity\\
\hline
\end{tabular}  
\end{center} 
\vskip 0.2cm

As above, $R\subset \R^n$ is a root system (irreducible and 
reduced), $W$ denotes the Weyl group
$<s_i, 1\le i\le n>,\,\, P$ is the weight lattice. 
\vskip 0.3cm

We omit the general definition of DAHA (it will be
provided later for $A_1$); see \textsuperscript{\cite{C101}}. 
 The following will be sufficient.
For $T_w$ as above,
$$\HH\!=\!\lan X_b,T_w,Y_b,q,t\ran,\, b\in P,w\in W,\ \,
\R\ni t^{\pm1/2},\, q\!=\!\exp(-\frac{1}{a}), a>0,$$
where the ring of coefficients is $\R$. More formally,
it is defined over  $\Z[q^{\pm\frac{1}{m}},t^{\pm\frac{1}{m}}]$
for proper $m$.
\vskip 0.2cm
 
\begin{definition}
The \textcolor{blue}{\sf \mathversion{normal}Shapovalov
anti-involution $\varkappa$ of $\HH$ for $Y$}
is such that $T_w^{\varkappa}\!=\!T_{w^{-1}}$ and 
the following ``PBW property" holds: for any $H\in \HH$,\,
the decomposition $H=\!\sum c_{awb}\! Y_a^{\varkappa} 
T_{w} Y_b$ exists and 
is unique.\sq
\end{definition}

{\sf An example}. Let $\varkappa: X_b\!\leftrightarrow\! 
Y_b^{\!-\!1}, T_{w}\!\to\! T_{w^{-1}} (w\!\in\! W).$
{\em All\,} Macdonald conjectures follow from its mere 
existence (without using the shift operator,
and practically without any calculations); 
see \textsuperscript{\cite{C6}}. 

\begin{definition}
The \textcolor{blue}{\sf coinvariant} is {\mathversion{normal}
$\{H\}_{\varkappa}^\varrho$} $\,\equal
\,\sum c_{awb}\, \varrho(Y_a)
\varrho(T_{w}) \varrho(Y_b)$,

\noindent
where ``PBW" is used, 
$\varrho$ is a linear map $\R[T_w,Y_b, w\in W,b\in P]
\to \R$ such that
$\rho:\R[Y^{\pm 1}]\to \R$
is a (one-dimensional)
character and $\varrho(T_w)\!=\!\varrho(T_{w^{-1}})$. 
A variant is with $\C$ instead of $\R$. Then
$\{\varkappa(H)\}_{\varkappa}^\varrho= 
\{H\}_{\varkappa}^\varrho$ by construction and \ 
$\{A,B\}\equal \{A^{\varkappa}\, B\}_{\varkappa}^\varrho
\,=\,\{B,A\}.$ \sq
\end{definition}

\textcolor{red}{\sf General problem.\,} Find an 
{\em integral\,} 
(analytic)
formula for $\{H\}_{\varkappa}^\varrho$. It is 
well defined for any $q,t\in \R^*$ (or in $\C^*$) by construction,
but presenting this ``algebraic" functional ``analytically"
is important in DAHA theory (and the key in this paper).


We will stick to the {\em polynomial case\,} through this 
paper. Namely, 
$\varrho$ will be the one-dimensional character of affine
Hecke algebra 
$\h_{Y}$ generated by $T_w$ and $Y_b$, which sends
$T_i\! \mapsto\! t^{1/2}, 
Y_b\!\mapsto\! t^{(\rho,b)}$ for $i\ge 0, b\in P$. 
Here $\rho=\frac{1}{2}\sum_{\al>0}\al$. Generally, the
number of different parameters $t$ here equals the number
of different lengths $|\al|$ in $R$.
Then $\{A,B\}$ acts via \,$\x\!\times\! 
\x$ for the {\sf polynomial representation}\,
$\x=\R[X^{\pm 1}]=Ind_{\h_Y}^{\HH} (\varrho)$.
\vskip 0.1cm

Generalizing the above definition, 
\textcolor{blue}{\sf\mathversion{normal}level-one
anti-involutions $\varkappa$\ } are such that 
dim$\HH/(\j\!+\!\j^\varkappa)\!=\!1$ for
$\x\!=\!\HH/\j$, $\j\!=\!\{H\mid H(1)=0\}$, $1\in \x$.
The Shapovalov ones are obviously level-one. Then
$\{H\}_{\varkappa}^\varrho$ is defined as the image of 
$H$ in $\HH/(\j\!+\!\j^\varkappa)$.

{\sf An example\,.} Let $\ast: g\mapsto g^{-1}$ for
$g=X_a,Y_b,T_w,q,t$. It is level-one for
\textcolor{red}{\it generic\,} $q,t$\,, but obviously
not a Shapovalov anti-involution with respect to  $Y$. One can 
prove \textsuperscript{\cite{C101}} that there exists 
the corresponding unique inner product in $\x$ for generic 
$q,t$ (not for all $t$ if $q$ is generic). 


\section{\textcolor{blue}{\sc Rational DAHA}}\label{SEC:RAT}
For rational DAHA, the counterpart of
$\ast$ above (serving the ``standard" inner product
in $\x$) is \textcolor{red}{\it\,not\,} level-one. 
The rational DAHA is:
 
$\HH''\equal\langle x,y,s\rangle/\,\{\,
[y,x]\!=\!\frac{1}{2}\!+\!k s,\ s^2\!=\!1,\
sxs\!=\!-x,\ sys\!=\!-y\,\}$.
\vskip 0.2cm

Accordingly the \textcolor{blue}{\sf
polynomial representation} $\x$ becomes $\R[x]$ with
the following action of $\HH''$:

\centerline{
$\!s(x)=-x,\ \,x=$ multiplication by $x,$\ \, $y\mapsto D/2,$}

\centerline{
where \ $D=\frac{d}{dx}+\frac{k}{x}(1-s)$ (the Dunkl operator). }
\vskip 0.2cm

%
%

Then the anti-involution $x^*\!=\!x, y^*\!=\!-y, s^*\!=\!s$ 
formally serves the inner product
$\int f(x)g(x)|x|^{2k}$, but it diverges
at $\infty$. Algebraically,
$\R[x]$ has \textcolor{red}{\it\,no\,} $*$-form for
$k\not\in -1/2-\Z_+$. Indeed, for $p\in \Z_+$
$\{1,y(x^{p})\}\!=\!0\!=\!\{1,c_{p}x^{p-1}\}$, where 
$c_{2p}\!=\!p,\, c_{2p+1}$
$=p+1/2+k$ (direct from the Dunkl operator).
Hence, $\{1,x^{p}\}=0$ ($\forall p$)
for non-singular $k$ and $\{\,,\,\}=0$. 

\vspace*{0.2cm}

To fix this problem, let us 
replace $y$ by $y+x$; then $*$ becomes Shapovalov
for such new $y$ (the definition depends on the choice
of $y$).
Indeed, the decomposition
$h=\sum c_{a\de b}((y+x)^*)^a s^\de (y+x)^b$ 
exists and is unique ($\de=0,1$) for any $h$.
Defining the {\sf coinvariant\,} by
 $\{h\}\equal\sum_{\de=0,1} 
c_{o\de o},\ \{f,g\}\equal\{f^*g\}$, it acts 
through $\R[x]e^{-x^2}\!\times \R[x]e^{-x^2}$ due to 
$(y+x)e^{-x^2}=0$ for the natural action of $\HH''$ on
$e^{-x^2}$. Indeed,
$\R[x]e^{-x^2}$ can be identified with
$\HH''/(\HH''(y+x),\HH''(s-1))$.
\vskip 0.3cm

Explicitly, let  $p\!=\!\frac{a+b}{2}$ for $a,b\in Z_+$. Then
a direct PBW calculation readily gives that
$\, \{x^a ,x^b \}\!=\!
(\frac{1}{2})^p(\frac{1}{2}+k)\cdots (\frac{1}{2}+k+p-1).$
Analytically, we ensured the convergence of 
$\int_{\R}fg|x|^{2k}$ via the
multiplication of $f$ and $g$
by $e^{-x^2}$; let us provide the exact analysis. 

\vskip 0.3cm

\textcolor{blue}{\it The integral presentation\,}
 for this form is: 
$$\{f\,,\,g\}=
\frac{1}{i}\int_{-\ep+i\R}(fge^{-2x^2}\,(x^2)^k)dx/
(\cos(\pi k)C),$$
where $C=\Gamma(k+1/2)\,2^{k+1/2}, 
\ \textcolor{red}{\it \forall k\in \C},\ \ep>0.$

\vskip 0.2cm

For \textcolor{blue}{\it real $k>-\frac{1}{2}\,$}, one can 
simply do the following:

$$\{f\,,\,g\}
=\frac{1}{i C}\int_{i\R} fge^{-2x^2}|x|^{2k}dx.$$
Note using $|x|$ here, which is not natural algebraically;
one can take here $x^{2k}$  instead using
the technique of hyperspinors (see below).
\vskip 0.2cm

Let \textcolor{blue}{\it $k\!=\!-\frac{1}{2}\!-\!m\,
(m\!\in\! \Z_+)$}. 
Then we replace
$\int_{-\ep\!+i\R}\rightsquigarrow 
\frac{1}{2}(\int_{-\ep\!+i\R}\!+\!
\int_{\ep\!+i\R})$ and
\ $\{f\,,\,g\}$ becomes
$\hbox{const\, Res}_0\,(fge^{-2x^2} x^{-2m-1}dx).$
The radical of this form is non-zero.
It is $(x^{2m+1}e^{-x^2})$, which is a 
\textcolor{red}{\it unitary\,} $\HH''$-module  
with respect to the form 
$\frac{1}{i}\int_{i\R} fge^{-2x^2}|x|^{-2m-1}dx$
restricted to this module (the convergence at $x=0$
is granted). The $*$-form of the quotient $\R[x]/(x^{2m+1})$\ 
is non-positive.
See \textsuperscript{\cite{CMa}} for some details.


\section{\textcolor{blue}{\sc General DAHA ($A_1$)}
\textsuperscript{\cite{C101}}}

\textcolor{blue}
{\sf The $q$-Extended elliptic braid group.} It is 
\ $\mathcal{B}_q\equal$ 

$\langle T,X,Y,q^{1/4}\rangle/$ 
$\{$\textcolor{blue}{\it\mathversion{normal}
$TXT\!=\!X^{-1},
TY^{\!-1}T\!=\!Y,
Y^{\!-1}X^{\!-1}YXT^2\!=\!q^{-\frac{1}{2}}$}
$\}.$
\vskip 0.2cm

\textcolor{blue}
{\sf Elliptic braid group.} It is  
$\mathcal{B}_1\equal\mathcal{B}_{q=1}=$
$\mathcal{B}_1=\pi_1^{\hbox{\eightrm orb}}(\{E\setminus 0\}/
\mathbf{S}_2),$ 
where $E$ is an elliptic curve (a $2$-dimensional torus).
We will provide below the geometric-topological interpretation of
the relations in $\mathcal{B}_1$.

\vskip 0.2cm

\textcolor{blue}
{\sf DAHA\,}  is defined as follows:\ \ {\mathversion{normal}
$\HH\equal\R[\mathcal{B}_q]/((T-t^{1/2})(T+t^{-1/2})),$} 
where $q=\exp(-1/a), a>0,\,  t=q^k.$ Here
$k\!\in\! \R$ for the positivity questions, but we
will need $k\in \C$ when doing analytic continuations. 
\vskip 0.3cm

If $t^{\frac{1}{2}}=1,$ then $T^2=1$ and we will replace
$T$ by $s$. In this case, $\HH$ becomes  
the \textcolor{blue}{ Weyl algebra} extended by  
$\mathbf{S}_2\, $. I.e. the relations are:
{
\begin{align*}
&sXs\!=\!X^{-1}, sYs\!=\!Y^{-1},
Y^{-1}X^{-1}YX\!=\!q^{-1/2},\, s^2=1.
\end{align*}
}

\noindent
Thus DAHA unites Weyl algebras 
with the Hecke ones. The Heisenberg and Weyl algebras 
(also called non-commutative tori) are the main 
tools in quantization of symplectic varieties. 
So DAHA can be expected to serve    
``refined quantization" (with extra parameters)
of  varieties with  global or local (in tangent spaces)
$W$-structures for Weyl groups $W$.  
\vskip 0.2cm


\begin{figure*}[htbp]
\vskip -0.5in
\hskip -0.5in
\includegraphics[scale=0.3]{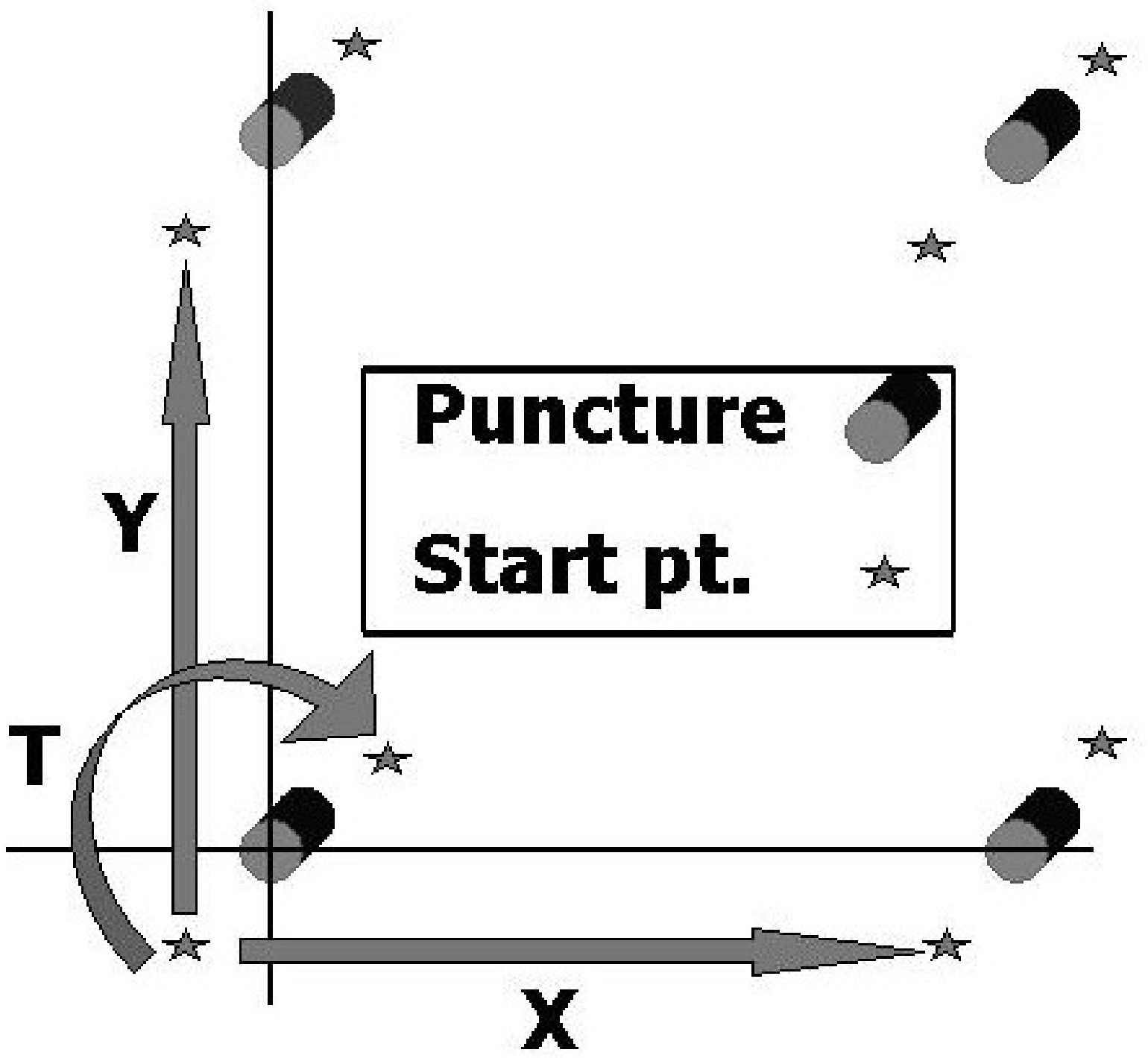}
\vskip -0.05in

\textcolor{blue}
{\it\mathversion{normal}
Generators  of $\b_1$ and 
relation $Y^{\!-\!1}X^{\!-\!1}YXT^2\!\! =\!\! 1$}
\vskip -0.0in

\hskip -0.4in
\includegraphics[scale=0.3]{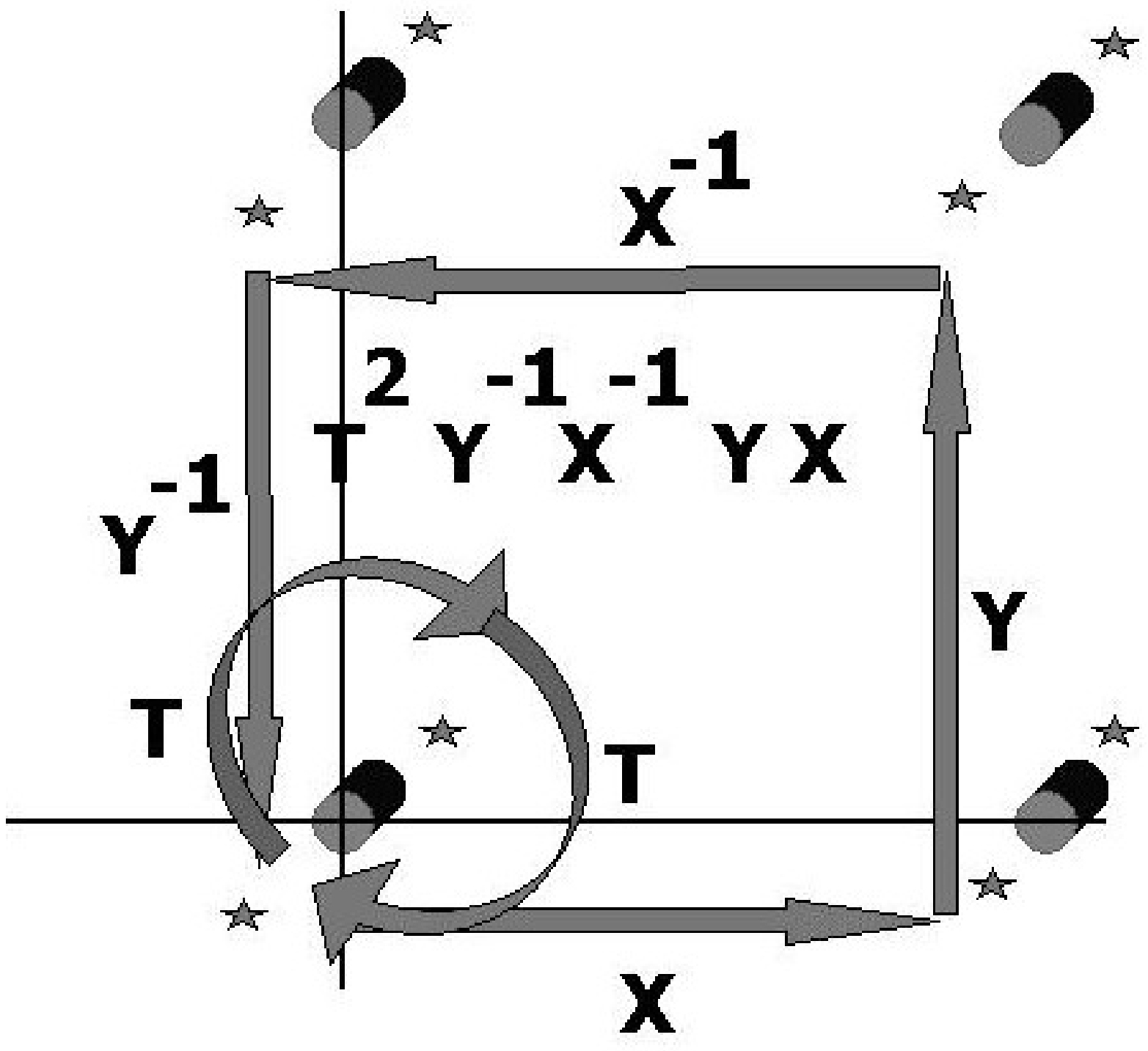}


\end{figure*}


The whole $PSL_2(\mathbf{ Z})$ acts projectively 
in $\b_q$ and $\HH$:
\begin{align*}
&\binom{11}{ 01}\sim \tau_+: 
Y\mapsto q^{-1/4}XY,\ X\mapsto X,\ T\mapsto T,\\ 
&\binom{10}{ 11}\sim \tau_-: X\mapsto q^{1/4}YX,\,\ \, 
Y\mapsto Y,\ \, 
T\mapsto T. 
\end{align*}

They are directly from topology. The key for us is
a pure algebraic fact that 
$\tau_+$  is the conjugation by $q^{x^2},$ where $X=q^x$;
use $\x$ below to see this.\ 
DAHA FT is for\ \, 
$\tau_+^{-1}\tau_-\tau_+^{-1}=\si^{-1}=$
$\tau_-\tau_+^{-1}\tau_-.$

\vskip 0.2cm
More exactly, 
the \textcolor{blue} {\it operator Fourier transform} is the 
DAHA automorphism sending:\, $q^{1/2}\mapsto q^{1/2},
t^{1/2}\mapsto t^{1/2}$,

$$Y\mapsto X^{-1},\ X\mapsto T Y^{-1}T^{-1},\ T\mapsto T;$$

\noindent
\textcolor{red}
{\it topologically, it is essentially
the transposition of the periods of $E$,\,}
though it is not an involution; it corresponds to the 
matrix {\tiny 
$\begin{pmatrix}0 & -1\\ 1 &0\end{pmatrix}$} representing 
$\si^{-1}$.

\comment{
\textcolor{blue} {\sf Rational DAHA} $\HH''$  
($[y,x]=1/2+ks$,\ \ldots)\ :

\noindent
$X\!=\!e^{\sqrt{\hbar} x}, Y\!=\!e^{-\sqrt{\hbar} y}, 
q\!=\!e^\hbar, t\!=\!q^k, \hbar\to 0, T\!\to\! s$.
}

\vskip 0.2cm
\textcolor{blue}
{\sf Polynomial representation.} It was defined above as
$Ind_{\h_Y}^{\HH}(\varrho).$ It is in the space
{\mathversion{normal}
$\x$} of {Laurent polynomials}  of $X=q^x.$ The action
is: 
\textcolor{blue}
{\it\begin{align*}
&T\mapsto \,t^{1/2}s\,+\,\, \frac{t^{1/2}-t^{-1/2}}
{q^{2x}-1}(s-1),\,\, Y\mapsto \pi T,\\ 
&\hbox{where\, \ }
\pi=sp, sf(x)=f(-x), s(X)=X^{-1},\\ 
&pf(x)=f(x+1/2),\ p(X)=q^{1/2}X, \ t=q^k.
\end{align*}}

\noindent
Here $Y$ becomes the \textcolor{blue}
{\sf difference Dunkl Operator}; $X$ acts by the multiplication.
\vskip 0.2cm

\noindent
The standard AHA stuff (Bernstein's Lemma)
gives that $Y+Y^{-1}$ preserves $\x_{sym}\equal$ 
$\{$ symmetric (even) Laurent polynomials$\}$;
the difference operator
$Y+Y^{-1}\mid_{sym}$ is sometimes called
the \textcolor{blue}{\sf $q,t$-radial part}.


\vskip 0.2cm
\textcolor{blue}
{\sf Basic inner products.}\textsuperscript{\cite{C101}}
Note that we {\it do not\,} conjugate $q,t$ below
(a simplest way to supply $\x$ with an inner product).  
\vskip 0.2cm

For $X\!=\!q^x,q\!=\!\exp(-\frac{1}{a})$ and
the Macdonald truncated $\theta$-function 
$$\mu(x)=\prod_{i=0}^{\infty}\frac{(1-q^{i+2x})(1-q^{i+1-2x})}
{(1-q^{i+k+2x})(1-q^{i+k+1-2x})},\,
\hbox{\, we set:}$$

\textcolor{blue}
{\it$\lan f,g \ran_{{1/4}}\!\equal\! 
\frac{1}{2\pi a i}\int_{1/4\!+\!P}\,f(x)\,T(g)(x) \mu(x) dx,\,$
where $P=[-\pi i a,\pi i a].$}
\vskip 0.3cm

\begin{theorem}\label{thm5.1}
For $k\!>\!-\frac{1}{2}$ (generally, $\Re k\!>\!
-\frac{1}{2}$), $\lan\!f,g\ran_{{1/4}}\!=$
$(fT(g)\mu)
_{\hbox{\sc\tiny ct}}$. The later inner product in $\x$ 
serves for any $k$ the anti-involution\, $\diamond$\,:
$ T^\diamond=T,\, Y^\diamond=Y,\, X^\diamond=T^{-1}XT$. The 
inner product  $\lan\!f,g\ran_{{1/4}}$ does it only for
$k\!>\! -1/2$, where it is {\it positive definite\,} 
in $\x=\R[X^{\pm 1}]$; however it remains symmetric 
for any $k$.  
\end{theorem}  
\vskip 0.2cm 

Proof. The coincidence of two formulas
for $\Re k>-1/2$ and the fact that $\diamond$ serves  
$(fT(g)\mu)_{\hbox{\sc\tiny ct}}$ are from 
\textsuperscript{\cite{C101}} (for any root systems).
The positivity is straightforward
via the norm-formulas for $E$-polynomials; let us 
provide a directly proof 
using that $\pi(\mu)=\mu(1/2-x)=\mu$.

\vskip 0.2cm

\noindent \textcolor{blue}{\sf (a)} The $E$-polynomials are
defined as follows:
$Y(E_n)=q^{-n_\#}E_n$,  where 
$n_\#= \frac{n-k}{2} \hbox{\, as\, } n\le 0,\ 
n_\#=\frac{n+k}{2} \hbox{\, as\, } n>0$,
$\ E_n=X^n+ (lower\  terms).$ Here $X_m$ is lower than
$X_n$ if either $|m|<|n|$ or $m=-n>0$.  

\smallskip
\noindent \textcolor{blue}{\sf (b)} Then
$\lan E_n,E_m \ran_{_{1/4}}=C_n\delta_{nm}$ for
some constant $C_n$ due to $Y^\diamond=Y$.

\smallskip
\noindent \textcolor{blue}{\sf (c)}
One has: $C_n=q^{-n_\#}\frac{1}{i}
\int_{1/4+P}\,E_n \overline {E_n} \mu(x) dx\,>\,0,$
\noindent
since $\pi (x)=\overline{x}$ (the latter is
complex conjugation) and $\mu(x)>0$ at $1/4+P$; 
use that 
$T(E_n)=\pi Y(E_n)=q^{-n_\#}\pi(E_n)=q^{-n_\#}\overline{E_n}.$
\sq
\vskip 0.2cm

\textcolor{blue}{\sf Imaginary Integration.} For
$\Re k\!>\!-\frac{1}{2}$, $f,g\!\in\!\x$, we set 
 
\textcolor{blue}
{\it$\lan f,g\ran_{\!_{1/4}}^{\ga,\infty}\equal\frac{1}{i}
\int_{\frac{1}{4}+i\R}\,fT(g)q^{-x^2} \mu(x) dx$}
\vskip 0.1cm

{\hfill $=\frac{1}{2i\sqrt{\pi a}}
\int_{\frac{1}{4}+P}\,fT(g)\sum_{j=-\infty}^{\infty}q^{j^2/4+jx}
\mu(x) dx$,}

\noindent
where we use that $f,T(g)$ and $\mu$ are $P$-periodic and  
that $q^{-x^2}$ is such with a multiplier. Then we 
employ the functional equation for the theta-function;
see Section 2.2.2 from \textsuperscript{\cite{C101}}.
For such $k$, this inner product
is symmetric and positive (as $k>-1/2$); it serves 
the anti-involution 
$$\varkappa :\ \ T^{\varkappa}=T,\, X^{\varkappa}=T^{-1}XT,\, 
Y^{\varkappa}=q^{-1/4}XY.$$
The latter involution is $\diamond$ above conjugated by
$\tau_+$, which reflects the multiplication the integrand
by $q^{-x^2}$. So the relation to $\varkappa$ 
(only for $\Re k>-1/2$) and 
the positivity follow from Theorem \ref{thm5.1}. 


%
%
%
%

\vskip -0.5cm

\section{\textcolor{blue}{\sc Analytic continuation}}
\vskip -0.1cm

The ingredients are as follows: the \textcolor{blue} 
{\it{Shapovalov}} $\,\varkappa\,$ above (for $Y$) and
the standard coinvariant $\varrho$ (serving $\x$). Recall
that 
$$\varrho\,(\sum_{a,b\in \Z}^{\ep=0,1}\, c_{a\ep b} 
(Y^\varkappa)^a\, 
T^\ep\, Y^b)\equal\sum c_{a\ep b} t^{\frac{a+\ep+b}{2}}$$
and the corresponding form is

\smallskip
\centerline {\mathversion{normal}
 $\{A,B\}_\varkappa^\varrho\,\equal\,\varrho(A^\varkappa B)=
\{B,A\}_\varkappa^\varrho$ in $\HH\ni A,B$.}

\smallskip
\noindent
The latter acts via \,$\x\! \times\! \x$, 
$\x\!=\!\R[X^{\pm 1}],$ 
and satisfies the normalization $\{1,1\}\!=\!1$ by construction. 
 \textcolor{red}{\it This form is 
regular (analytic) for all $k\in\C$\,}.
\vskip 0.3cm

\begin{theorem}
For $\Re k >-1/2,$ one has:
$G(k)\{f,g\}_\varkappa^\varrho=\,
\lan\!f,g\ran_{\!_{1/4}}^{\ga,\infty},$
where
$G(k)=
\sqrt{\pi a}\prod_{j=1}^\infty\,\frac{1-q^{k+j}}{1-q^{2k+j}}$\,
(the latter is from \textsuperscript{\cite{C101}},Theorem 2.2.1). 
\end{theorem}

{\it Proof}.  Let
{\mathversion{normal}
\textcolor{blue}{\it$\Phi_\ep^k(f,g)$}$\equal$
$\frac{1}{i}\int_{\ep+i\R}\,fT(g)q^{-x^2} \mu(x) dx$}
for the path $\c\equal\{\ep+i\R\}$. For such a path, 
{\it bad } (singular) $k$ are 
$\{\,2\c-1-\Z_+,\, -2\c-\Z_+\,\}$ (when poles
of $\mu$ belong to $\c$);
so $\{\Re k>-1/2\}$ are all \textcolor{red}{\it good } 
as $\ep=1/4$. Then we use that the theorem holds for
$\Re k\!>\!\!>0$.
\sq

\vspace*{0.3cm}

\textcolor{blue}{\sf The case $\ep\!=\!0$.}  
Then $\Phi_0^k(f,g)$ coincides
with $G(k)\{f,g\}_\varkappa^\varrho$ only for $\Re k>0$. 
\textcolor{red}{\it For any $k$\,}, this form 
is symmetric and
its anti-involution sends
$T\!\!\mapsto\!T,$ $X\!\mapsto\!X^{\!\varkappa}\!=\!\!T^{-1}XT$ 
(the image of $Y$ is {\it not\,} $Y^\varkappa$ if $\Re k<0$).
\vskip 0.1cm

\textcolor{blue}{\sf Comparing  
$\ep\!=\!0$ and $\ep\!=\!\frac{1}{4}$ for
$0\!>\!\Re k\!>\!-\frac{1}{2}.$}\ ({
This is actually the induction
step for the analytic continuation to any negative $\Re k$}).

\textcolor{red}{\it Let us assume that  
$F=fT(g)\in \R[X^{\pm 2}]$}. We mainly follow 
\textsuperscript{\cite{CO}}, Section 2.2.
By picking up the residues between 
the $\c$-paths at $0$ and at $1/4$, one obtains that
$\Phi_{\frac{1}{4}}^k\!=\!
\Phi_0^k\!+\!A(-\frac{k}{2})\mu^\bullet
(-\frac{k}{2}) F(-\frac{k}{2})
\equal \hat{\Phi}^k,$
where
$A(\tilde{k})=\sqrt{\pi a}\,
\sum_{m=-\infty}^{\infty}q^{m^2+2m\tilde{k}}$\ 
(the contribution of $q^{-x^2}$), and

\vskip 0.1cm
$F(-\frac{k}{2})=fT(g)(x\mapsto\! -\frac{k}{2}),$\ \ 
$\mu^\bullet\,(-\frac{k}{2})=\prod_{j=0}^{\infty}\, \frac
{(1-q^{k+j+1})(1-q^{-k+j})}
{(1-q^{1+j})(1-q^{2k+j+1})}=$
\vskip 0.1cm

\noindent
{\small $\bigl((1-q^{2x+k})\mu(x)\bigr)(x\mapsto -k/2).$}
\vskip 0.2cm

Importantly,
$\hat{\Phi}^k$
is meromorphic for \textcolor{red}{\it $\Re k\!>\!-1$} 
(i.e. beyond $-\frac{1}{2}$ for $\Phi_{\frac{1}{4}}^k$); so 
it coincides with $G(k)\{f,g\}_\varkappa^\varrho$ there. 
We note that $\hat{\Phi}^k$ is symmetric 
\vskip 0.1cm

\noindent
for {\it any\,} $k$. Indeed:\ \ 
$fT(g)(-\frac{k}{2})\,=\,
t^{1/2}f g(-\frac{k}{2})\,=\,T(f)g(-\frac{k}{2})\ $
due to
\smallskip

\noindent
$T=\frac{q^{2x+k/2}\!-q^{-k/2}}{q^{2x}-1}\,s-
\frac{q^{k/2}\!-q^{-k/2}}{q^{2x}-1}$
and $\ (q^{2x+k/2}-q^{-k/2})(x\mapsto -k/2)=0$.

\vskip 0.3cm

\begin{maintheorem} 
For  $F=fT(g)\in \R[X^{\pm 2}]$
and for any $\Re k<0$\,:
\vskip 0.1cm

\noindent
$G(k)\{f,g\}_\varkappa^\varrho=\Phi_0^k$
$+\,\mu^\bullet(-k/2)\sum_{\tilde{k}\in \tilde{K}}\,
A(\tilde{k})\,\left(
\mu^\bullet(\tilde{k})/\mu^\bullet(-\frac{k}{2})\right)
\,F(\tilde{k})$, for 

\noindent
$\tilde{K}=\{n_\#,|n|\le m\}=\{-k/2\}\cup\{\pm\frac{k+j}{2},\,
1\!\le\!j\le m\},\ \,m\equal[\Re(-k)],$
\vskip 0.1cm

\noindent
where $[\cdot]\!=$integer part, and\ 
$\mu^\bullet(\pm\frac{k+j}{2})/\mu^\bullet(-\frac{k}{2})=
t^{-j_\pm}
\prod_{i=1}^{j_\pm}\frac{1\!-\!t^2q^i}{1\!-\!q^i}$  
for $j_+\!=\!j\!-\!1, j_-\!=\!j$.
Generally, if $F\in\R[X^{\pm 1}]$ (not in $\R[X^{\pm 2}]$
as above), the poles of 
$\mu$ are given by the relations
$q^{-\frac{1}{2}} X \in$ 
\textcolor{red}{\it\mathversion{normal} $\pm\,$} 
$q^{\Z_+/2}\,t^{\frac{1}{2}}
\ni X^{-1}$, 
and the summation must be ``doubled" accordingly.   
\sq
\end{maintheorem}

Here we count ``jumps" through the walls 
$\Re k\!=\!-j\!\in\! -\Z_+$. The duplication of the summation
for $F\in\R[X^{\pm 1}]$ corresponds to the passage from 
affine Weyl group $\tW$ to its extension $\hW$ by $\Pi=\Z_2$
in the $p$-adic limit $q\to 0, X\mapsto Y$ 
(discussed below).

\vskip 0.2cm
Importantly, here and for any root systems 
only the total sum is an $\HH$-invariant form.
The partial sums with respect to the dimensions
of the integration domains are symmetric and even 
$\h_X$-invariant, but they are {\em not\,} $\HH$-invariant.

\vspace*{0.2cm}
\noindent

\begin{corollary} The form
 $\{f,g\}_\varkappa^\varrho$ is degenerate exactly at
the poles of $G(k): k\!=\!-\frac{1}{2}-m, m\!\in\! \Z_+$. For such
$k$, the quotient of $\x$ by its radical is a direct
sum of $2$ irreducible
$\HH$-modules of dim$=2m+1$ (``perfect" in the
terminology from \textsuperscript{\cite{C101}}), transposed
by the map $X\mapsto -X$.
\sq
\end{corollary}
\vskip 0.1cm

\textcolor{red}
{
The radical here is the ideal $(E_{2m+1})$, which is
a {\it unitary\,}  $\HH$-module with respect 
to $\Phi_0^k$,\,} matching  the 
analogous fact for the rational DAHA  $\HH''$ observed above. 
The rational limit  is as follows:
\noindent
$q\!=\!e^{\hbar}, t\!=\!q^k,
Y\!=\!e^{-\sqrt{\hbar} \bar{y}}, X\!=\!e^{\sqrt{\hbar} \bar{x}},
\, \hbar\!\to\! 0$
($\bar{x}$ is {\it not\,} $x$ from $X=q^x$).
One has $q^{x^2}= e^{\bar{x}^2}$, and $\mu(x)$ becomes
essentially  $\bar{x}^{2k}$ under this limit. 
So the space
Funct$(\tilde{K})$ (upon the restriction to $F\in \R[X^{\pm 2}]$, 
which does not influence the rational limit) 
directly maps to the $\HH''$-module 
$\R[\bar{x}]/(\bar{x}^{2m+1})$ (in terms of $\bar{x},\bar{y}$).
\begin{theorem}
For any root system $R\subset \R^n$, 
a Shapovalov or level-one anti-involution $\varkappa$,
and for the coinvariant $\varrho$ serving the polynomial 
representation, 
the corresponding DAHA-invariant form can be represented as
a finite sum of integrals over translations of 
$\imath\mathbb{A}$ for proper subspaces $\mathbb{A}\subset \R^n$, 
called double-affine $q,t$-residual subtori, starting with
the full imaginary integration. \sq  
\end{theorem} 

\vskip 0.3cm
The program is to  $(a)$\  find (explicitly) 
these subtori for any root systems,  
$(b)$\  calculate the corresponding 
$C$-coefficients (the $q$-deformation
of the AHA Plancherel measure),
and finally $(c)$\ perform the $p$-adic 
limit ($q\to 0$), which steps are
technically non-trivial (even in type $A$). 


\section{\textcolor{blue}{\sc P--adic limit}}

For the AHA $\h$ of type $A_1$, we set
$s=s_1, \om=\om_1, \pi=s\om$. Let 
$$\psi_n\!\equal\! t^{-\frac{|n|}{2}}T_{n\om}\p_+,\, 
\p_+\!=\!(1+t^{1/2}T)/(1\!+\!t) \for n\in \Z.$$
One can naturally consider them 
as  polynomials in terms of $Y\equal T_{\om}=\pi T$; 
then they become the {\it Matsumoto spherical functions}. 
This identification is based on the analysis by Opdam 
\textsuperscript{\cite{O12}} and the 
author \textsuperscript{\cite{CO}}; see also 
\textsuperscript{\cite{C101}}(Section 2.11.2) and
\textsuperscript{\cite{Ion,CMa}}.
Accordingly, 
the {\it Satake-Macdonald $p$-adic spherical functions} 
become $\,\p_+\psi_n (n\ge 0)$. 
 
\vskip 0.2cm

\begin{theorem} 
For $n\in \Z$,\, the polynomials $E_n(X)/E_n(t^{-\frac{1}{2}})$
become $\psi_n$ as $q\!\to\!0$ upon the following
substitution:

\textcolor{red}{\it\mathversion{normal}
$$
f(X)\mapsto f(X)'
\!\equal\!
f(X\mapsto X'=Y,\,  t\mapsto t'=\frac{1}{t}).$$}

\noindent
Let $\mu_0\!=\!\mu(q\to 0)\!=\!\frac{1-X}{1-tX}$,
$\{f,g\}_0\!=\!(f T(g)\mu_0)_{\hbox{\tiny\sc ct}}$. Then
for $\lan T_{\hw}\ran=\de_{id,\hw}$ and the standard
anti-involution $T_{\hw}^\star=
T_{\hw^{-1}}$ in $\h$, one has:

$$\{f,g\}_0(t\mapsto t')=(t^{1/2}\!+\!t^{-1/2})
\lan (f'\p_+) (g'\p_+)^\star\ran \for f,g\!\in\! \x,$$

\noindent
which is actually the \textcolor{red}{\it nonsymmetric\,} 
AHA Plancherel formula for the $p$-adic
Fourier transform. Here $t',f',g'$ are as above.\sq
\end{theorem}


The corresponding version of the Main Theorem 
(compatible with the $p$-adic limit) is as follows.
The Gaussian collapses and we must omit it and use the
integration over the period instead
of the imaginary integration. We continue
using the notations $j_{\pm}\!=\!
\{j-1,j\}, t=q^k$.

\begin{theorem}\label{NO-Gauss}
For $q=e^{-\frac{1}{a}}$,\, $M\!\in\! \N/2,\,\, 
F(x)=fT(g)(q^x)\!\in\R[q^{\pm 2x}]\,:$

\begin{align*}
(F\mu)_{\hbox{\tiny\sc ct}}
=&\,\frac{1}{2\pi a M \imath}
\int_{-\pi a M\imath}^{+\pi a M \imath}
F(x)\mu(x)\, dx\\
+\,\mu^\bullet(-\frac{k}{2})&\,\times
\Bigl(F(-\frac{k}{2})+ \sum_{j=1,\,\pm}^{[\Re(-k)]} 
F(\pm\frac{k+j}{2})
\ t^{-j_\pm}\prod_{i=1}^{j_{\pm}}
\frac{1\!-\!t^2q^i}{1\!-\!q^i}
\Bigr).
\end{align*}
Here $k$ is arbitrary. The left-hand side is entirely algebraic 
and meromorphic for any $\,k\,$ by construction.
 Namely \textsuperscript{\cite{C101}},
$(F\mu)_{\hbox{\tiny\sc ct}}=(\mu)_{\hbox{\tiny\sc ct}}
(F\mu^\circ)_{\hbox{\tiny\sc ct}}$,
\vskip -0.3cm
{\small
$$\mu^\circ\equal\mu(x)/(\mu)_{\hbox{\tiny\sc ct}}=
1+\frac{q^k-1}{1-q^{k+1}}(q^{2x}+q^{1-2x})+\cdots
$$
}

\vskip -0.2cm
\noindent
is a series in terms of {\small $(q^{2mx}+q^{m-2mx})$ } for
$m\ge 0$ with rational $q,t$-coefficients, 
which is essentially Ramanujan's ${}_1\!\Psi_1$-summation, and
{\small 
$$
(\mu)_{\hbox{\tiny\sc ct}}=\frac{(1-q^{k+1})^2(1-q^{k+2})^2\cdots }
{(1-q^{2k+1})(1-q^{2k+2})\cdots (1-q)(1-q^2)\cdots }\ .
$$
}
\end{theorem}
\vskip -0.5cm\sq

One can replacing the integral above
with the corresponding sum of the residues,
which is an interesting generalization of the classical
formula for the {\em reciprocal\,} of the theta-function
\textsuperscript{\cite{Car}}. Its extension to any root systems
requires {\it Jackson integrations\,}; see Section 3.5 from 
\textsuperscript{\cite{C101}}.

\begin{proposition}
For $\Re k<\!-m\!\in -\Z_+$ and\, $F(x)\!\in q^{-2m}\R[q^{+2x}]$, 
\begin{align*}
\frac{1}{\pi a\imath }\int_{-\pi a \imath/2}^{+\pi a \imath/2}
F(x)\mu(x)\, dx\ =\ &\\
\mu^\bullet(-\frac{k}{2})\,\times
\Bigl(-F(-\frac{k}{2})-\sum_{j=1}^{[\Re(-k)]} 
&\ F(-\frac{k\!+\!j}{2})
\ t^{-j}\prod_{i=1}^{j}
\frac{1\!-\!t^2q^i}{1\!-\!q^i}\\
+ \sum_{j=[\Re(-k)]+1}^{\infty} 
&\ F(\frac{k+j}{2})
\ t^{1-j}\prod_{i=1}^{j-1}
\frac{1\!-\!t^2q^i}{1\!-\!q^i}
\Bigr),\\
(F\mu)_{\hbox{\tiny\sc ct}}=
\mu^\bullet(-\frac{k}{2})\,\times
\Bigl(\sum_{j=1}^{\infty} 
&\ F(\frac{k+j}{2})
\ t^{1-j}\prod_{i=1}^{j-1}
\frac{1\!-\!t^2q^i}{1\!-\!q^i}
\Bigr).
\end{align*}
\end{proposition}
\vskip -0.9cm\sq


\vskip 0.2cm
Switching in (\ref{NO-Gauss})
to $X=q^x$ and making $a=\frac{1}{M}$ for $M\to \infty$
(then $q\!\to\! 0$), let 
$k\!=\!-ca$ for $c>0$. Then
$t\!=\!e^{-\frac{k}{a}}\!\to\! e^c$ and the formula
above under $\Re k\to 0_-$ becomes the Heckman-Opdam one;
recall that DAHA with $t\!>\!1$ is related to AHA 
from Section \ref{SEC:AHA} for $t'\!=\!\frac{1}{t}<1$.  
Here and for any root systems, 
\textcolor{red} {\it only AHA residual subtori\,} 
contribute for $M>\!>0$.

\section{\textcolor{blue}{\sc Conclusion}} 

Let us summarize the main elements
and steps of the construction we propose.
The ingredients are as follows.
\vskip 0.1cm
 
\noindent
(a)\, {\it Shapovalov
anti-involution} \ $\varkappa$\ of $\HH$ (with respect to 
the subalgebra $\y=\R[Y_b]$), i.e.\, such that
\,$\{\varkappa(Y_a) T_w Y_b\}$\, form a (PBW) basis of $\HH$; 

\noindent
(b)\, the corresponding {\it coinvariant\,}: $\varrho:\HH\to \R$
satisfying $\varrho(\varkappa(H))=\varrho(H)$ (for 
any character of $\y$ and $\varrho$ 
on $\H$ s.t. $\varrho(T_w-T_{w^{-1}})\!=\!0$);
 
\noindent
(c) the corresponding \,{\it Shapovalov form\,}  
$\{ f,g\}\equal$
$\varrho(A^\varkappa\, B)$ for $A,B\in \HH$,
satisfying $\{ 1,1\}=1$ and
analytic for any $k$.

\vskip 0.1cm
\textcolor{red}{\it The main problem\,} is to express 
$\{f,g\}_\varkappa^\varrho$ as a sum of integrals 
over the \textcolor{blue} {\it DAHA residual subtori}
for any (negative) $\Re k$. Then one can try
to generalize this formula  to arbitrary
DAHA anti-invo\-lutions (any ``levels") and any 
induced modules.

\vskip 0.1cm
\textcolor{blue}{\sf Hyperspinors}
\textsuperscript{\cite{Ch11,  CMa, ChO2, O2}}.
An important particular case of the program above is
a generalization of the integral formulas from
the spherical case to the whole regular representation of $AHA$.
The technique of {\it hyperspinors\,} is expected to
be useful here; they were called
$W$-{\it spinors} in prior works ($W$ stands for the Weyl group).

The $W$-spinors are 
simply collections $\{f_w,w\in W\}$ of elements
$f_w\in A$ with a natural action of $W$ on the
indices. If $A$ (an algebra or a sheaf of algebras) has its
own (inner) action of $W$ and $f_{w}=w^{-1}(f_{id})$, they
are called {\it principle spinors\,}.
Geometrically, hyperspinors are  
$\C W$-valued functions on any
manifolds, which is especially interesting
for those with an action of $W$.
The technigue of spinors can be seen as a direct generalization of
{\it supermathematics}, which is the case of the root system
$A_1$, from $W=\S_2$ to arbitrary Weyl groups.
\vskip 0.2cm

For instance,  Laurent polynomials with the
coefficients in the group algebra $\C W$ are considered 
instead of $\x$, the integration is defined
upon the projection $W\ni w\mapsto 1$ (a counterpart
of taking the even part of a super-function), 
and so on and so forth. No ``brand new" definitions are 
necessary here, but the theory quickly becomes involved. 

\vskip 0.2cm
The $W$-spinors  proved to be very useful for quite a few 
projects. One of the first instances was
the author's proof in \textsuperscript{\cite{Ch11}} of
the Cherednik-Matsuo theorem, an isomorphism between the
{\it AKZ\,} and {\it QMBP\,}. An entirely algebraic version of 
this argument was presented in \textsuperscript{\cite{O2}};
also see \textsuperscript{\cite{C101}}.
This proof included the
concept of the fundamental group for the configuration space
associated with $W$ or its affine analogs
{\it without fixing a starting point}, \`a la Grothendieck.
A certain system of cut-offs and the related
{\it complex hyperspinors\,}
can be used instead.
The corresponding representations of the braid 
group becomes a $1$-cocycle on $W$ (a much more algebraic object
then the usual monodromy). 
 
\vskip 0.2cm

A convincing application of the
technique of hyperspinors was the
theory of {\it non-symmetric} \,$q$-Whittaker functions.
The Dunkl operators in the theory of Whittaker functions
(which are non-symmetric as well as the corresponding
Toda operators) simply cannot be
defined without hyperspinors and the calculations with them
require quite a mature level of the corresponding technique.
See \textsuperscript{\cite{CMa}} and especially
\textsuperscript{\cite{ChO2}} (the case of arbitrary root
systems). The Harish-Chandra-type decomposition formula
for global {\it nonsymmetric\,} functions from 
\textsuperscript{\cite{C7}} (for $A_1$)
is another important application;
hyperspinors are essential here.

\vskip 0.2cm
By the way, $x^{2k}$ for complex $k$, which is one of
the key in the rational theory (see above), is a typical 
{\it complex
spinor\,}, i.e. a collection of two 
(independent) branches of this function in the upper and 
lower half-planes. To give another (related)
example, the Dunkl eigenvalue problem 
always has $|W|$ independent {\it spinor\,}
solutions; generally, only 
one of them is a {\it function}. In the case of $\HH''$ for
$A_1$ (above), both fundamental spinor solutions 
for singular $k=-1/2-m,\, m\in \Z_+$ are {\it functions}. See 
\textsuperscript{\cite{CMa}} for some details.
\vskip 0.2cm

A natural question is,  do we have hypersymmetric
physics theories for any Weyl groups $W$,
say ``$W$-hypersymmetric Yang-Mills theory"?

\vfil
\vskip 0.3cm
{\sf Jantzen filtration.}
\textcolor{red} It is generally  
a filtration of the polynomial representation
of  $\x$ in terms of {\it AHA modules\,}, not
DAHA modules, for $\Re k<0$. The top module 
is the quotient of $\x$ by the radical of the
sum of integral terms for the smallest residual
subtori (points in many cases). Then we restrict
the remaining sum to this radical and continue by
induction with respect to the dimension of 
the (remaining) subtori. 

Sometimes certain
sums for residual subtori of 
dimensions smaller than $n$ are DAHA-invariant; then 
$\x$ is reducible.  We expect that the reducibility of $\x$
always can be seen this way,  which
includes the degenerations of DAHA.
For $A_n$, the corresponding 
Jantzen filtration provides the whole decomposition
of $\x$ in terms of irreducible DAHA modules,
the so-called {\it Kasatani 
decomposition\,} \textsuperscript{\cite{En,ES}}.
Generally, the corresponding quotients can be DAHA-reducible.
\vskip 0.2cm
 
For instance,  the {\it bottom module\,} 
of the Jantzen filtration has the inner product
that is (the restriction of) the integration
over the whole $i\R^n$. This provides some {\it a priori\,}
way to analyze its signature (positivity), which is of obvious
interest. The bottom DAHA submodule of $\x$
was defined algebraically (without the Jantzen 
filtration) in \textsuperscript{\cite{C6}}.
Indeed, it appeared semisimple under certain technical 
restrictions. For $A_n$, this is related to the so-called 
{\it wheel conditions\,}. 
\vskip 0.2cm

Let us discuss a bit the \textcolor{blue} {\it rational case.}
The form $\{f,g\}_\varkappa^\varrho$ for $\HH''$
can be expected
to have a presentations in terms of integrals
over the $x$-domains with $\Re x$ in  
the boundary of a \textcolor{blue} {\it tube 
neighborhood} of the
\textcolor{blue} {\it resolution} of the cross 
$\prod_{\al\in R_+} (x,\al)=0$ over $\R$. The simplest
example is the integration 
over $\pm i\ep +\R$ for $A_1$. 
This resolution (presentation of the cross as a divisor with
normal crossings) is due to the author 
{\small (Publ. of RIMS, 1991}), 
de Concini - Procesi, and Beilinson - Ginzburg. 
\vskip 0.1cm

This can be used to study the {\it bottom module\,}
of the polynomial representation for 
{\it singular\,} $k_o=-\frac{s}{d_i}$,
assuming that it is well-defined and $\HH''$-invariant.
When $s\!=\!1$ (not for any $s$), it can be proved unitary
in some interesting cases; see 
Etingof {\it et al.} in the case of $A_n$ 
\textsuperscript{\cite{ES}}. The restriction
of the initial (full) integration over $\R^n$ 
provides a natural approach to this phenomenon. 
See Section \ref{SEC:RAT} in the case of $A_1$.

The DAHA-decomposition of the polynomial or other modules
is a natural application of the integral 
formulas for DAHA-invariant forms,
but we think that knowing such formulas is necessary
for the DAHA harmonic analysis even if 
the corresponding modules are irreducible.  

{\sf Acknowledgements.}
The author thanks RIMS, Kyoto university for the
invitation, and the participants of his course at UNC.
Many thanks to the referee for important remarks.


\vskip -3cm
\bibliographystyle{unsrt}


\end{document}

%% file: macro.tex
\renewcommand{\tilde}{\widetilde}
\renewcommand{\hat}{\widehat}

\newcommand{\Z}{{\mathbb Z}}

\newcommand{\N}{{\mathbb N}}
\newcommand{\C}{{\mathbb C}}
\newcommand{\R}{{\mathbb R}}

\def\HH{\mbox{${\mathcal H}$\kern-5.2pt${\mathcal H}$}}

\newtheorem{theorem}{Theorem}[section]
\newtheorem{maintheorem}[theorem]{Main Theorem}
\newtheorem{proposition}[theorem]{Proposition}
\newtheorem{definition}[theorem]{Definition}

\newtheorem{corollary}[theorem]{Corollary}

\newtheorem{theorem }{Theorem}[section]
\newtheorem{maintheorem }[theorem]{Main Theorem}
\newtheorem{proposition }[theorem]{Proposition}
\newtheorem{definition }[theorem]{Definition}
\newtheorem{lemma }[theorem]{Lemma}
\newtheorem{corollary }[theorem]{Corollary}
\newtheorem{notation }[theorem]{Notation}
\newtheorem{remark }[theorem]{Remark}
\newtheorem{example }[theorem]{Example}

\newtheorem{ maintheorem }[theorem]{Main Theorem}
\newtheorem{ theorem}{Theorem}[section]
\newtheorem{ proposition}[theorem]{Proposition}
\newtheorem{ definition}[theorem]{Definition}
\newtheorem{ lemma}[theorem]{Lemma}
\newtheorem{ corollary}[theorem]{Corollary}
\newtheorem{ notation}[theorem]{Notation}
\newtheorem{ remark}[theorem]{Remark}
\newtheorem{ example}[theorem]{Example}

\def\for{\  \hbox{ for } \ }

\def\and{\  \hbox{ and } \ }
\def\and{\  \hbox{ and } \ }

\def\equal{\stackrel{\,\mathbf{def}}{= \kern-3pt =}}

\def\la{\lambda}
\def\La{\Lambda}
\def\om{\omega}

\def\al{\alpha}

\def\ga{\gamma}
\def\ep{\epsilon}

\def\de{\delta}

\def\si{\sigma}

\def\tW{\widetilde W}

\def\hw{\widehat{w}}
\def\hW{\widehat{W}}

\def\S{\mathbf{S}}

\def\0{\mathbf{0}}
\def\H{\mathbf{H}}

\def\çF{\mathcal{F}}

\def\l{\mathcal{L}}

\def\p{\mathcal{P}}

\def\h{\mathcal{H}}
\def\c{\mathcal{C}}
\def\y{\mathcal{Y}}

\def\x{\mathcal{X}}

\def\j{\mathcal{J}}
\def\b{\mathcal{B}}

\def\lan{\langle}

\def\ran{\rangle}

\newcommand{\sq}{\phantom{1}\hfill$\qed$}

\def\HH{\mathfrak{H}}

\def\HH{\hbox{${\mathcal H}$\kern-5.2pt${\mathcal H}$}}
\def\HHH{\hbox{${\mathbb H}$\kern-4.2pt${\mathbb H}$}}

\font\smm=msbm10 at 12pt 
\def\symbol#1{\hbox{\smm #1}}
\def\lsmash{{\symbol n}}

\def\#{\sharp}

\font\eightrm=cmr8
